\documentclass[runningheads,a4paper]{llncs}
\usepackage{amsmath}                                         
\usepackage{amssymb}
\usepackage{graphicx}
\usepackage{url}

\usepackage[framed, hyperref, thmmarks, amsmath]{ntheorem}   
\theorembodyfont{\upshape}
\theoremstyle{plain}
\newtheorem{alg}{Algorithm}
\newtheorem{funct}{Function}
\begin{document}

\mainmatter

\newcommand{\add}[1]{}

\title{Numerical Software to Compute Newton Polytopes and Tropical Membership}  
\titlerunning{Newton Polytopes and Tropical Membership} 
\author{Taylor Brysiewicz\inst{1} }
\authorrunning{Brysiewicz}
\institute{
Texas A\&{M} University, United States of America\\
\email{tbrysiewicz@math.tamu.edu},\\ 
\texttt{http://www.math.tamu.edu/}$\sim$\texttt{tbrysiewicz/}
}
\maketitle
\begin{abstract}
We present our implementation of an algorithm which functions as a numerical oracle for the Newton polytope of a hypersurface in the {\bf Macaulay2} package {\bf NumericalNP.m2}. We propose a tropical membership test, relying on this algorithm, for higher codimension varieties based on ideas from Hept and Theobald. To showcase this software, we investigate the Newton polytope of both a hypersurface coming from algebraic vision and the L\"uroth invariant.
\end{abstract}
\section{Introduction}
Often hypersurfaces are presented as the image of a variety under some map. Determining the defining equation $f \in \mathbb{C}[x_1,\ldots,x_n]$ of such a hypersurface $\mathcal H \subseteq \mathbb{C}^n$ is computationally difficult and often infeasible using symbolic methods such as Gr\"obner bases. Moreover, many times the defining equation is so large that it is not human-readable and so one naturally desires a coarser description of the polynomial, such as the Newton polytope. The Newton polytope of $f$, or equivalently that of $\mathcal H$, is the convex hull of the exponent vectors appearing in the support of $f$ and provides a large amount of information about the hypersurface. Newton polytopes are necessary to compute the BKK bound on the number of solutions to a polynomial system \cite{BKK1} and can also provide topological information such as the Euler characteristic of the hypersurface \cite{NPAG}. Knowing $\textrm{New}(f)$ also reduces the computational difficulty of finding $f$ via interpolation: the size of the linear system one must solve is $|\textrm{New}(f) \cap \mathbb{Z}^n|$, which is usually much smaller than the na\"ive bound of ${n+d-1}\choose{d}$ where $d=\deg(f)$.

In 2012 Hauenstein and Sottile \cite{NP} proposed an algorithm we call the {HS-algorithm} (Algorithm \ref{HSAlgorithm}) and showed that this algorithm functions as a vertex oracle for linear programming on $\textrm{New}(f)$. This algorithm requires that the hypersurface is represented numerically by a witness set.  Because a witness set is the only requirement, the HS-algorithm applies to hypersurfaces which arise as images of maps. 

We observe that the HS-algorithm is stronger than a vertex oracle and so we introduce the notion of a numerical oracle which returns some information even when the linear program is not solved by a vertex.
This observation gives rise to an algorithm for determining membership in a tropical variety (Algorithm \ref{TropicalMembership}) based on ideas of Hept and Theobald \cite{HeptTheobald}.
Both the HS-algorithm and a prototype of the tropical membership algorithm have been implemented in the {\bf Macaulay2} \cite{M2} package {\bf  NumericalNP.m2} which uses the package {\bf Bertini.m2} \cite{B4M2}  to call {\bf Bertini} \cite{Bertini} to perfom numerical path tracking. 

Section \ref{Section:Theory} contains background on polytopes, numerical algebraic geometry, and tropical geometry. A description of both the HS-algorithm and the tropical membership algorithm along with bounds on the convergence rates involved in the tropical membership algorithm can be found in Section \ref{Sec:Algorithms}. Section \ref{Section:Functionality} outlines the main user functions in {\bf NumericalNP.m2} and Section \ref{Section:Applications} advertises the stength of the software on much larger examples.
\section{Underlying Theory}
\label{Section:Theory}
\subsection{Polytopes}
A \emph{polytope} $P\subseteq \mathbb{R}^n$ is the convex hull of finitely many points $V \subseteq \mathbb{R}^n$. Equivalently, $P$ is the bounded intersection of finitely many halfspaces. The former presentation is a \emph{$V$-representation} of $P$ while the later is an \emph{$H$-representation} of $P$. 
Given $\omega \in \mathbb{R}^n$ the set ${P_\omega}:=\{x \in P \mid \langle x, \omega \rangle \textrm{ is maximized}\}$ is called the \emph{face} of $P$ \emph{exposed} by $\omega$ and the function $h_P(\omega)=\underset{x \in P}{\max} \langle x,\omega \rangle$ 
is the \emph{support function} of $P$. We define a \emph{numerical oracle} to be the function
$$ {\mathcal O}_P: \mathbb{R}^n \to  \mathbb{N}^n \cup \{\texttt{EEP}\} $$
$$\omega \mapsto \begin{array}{cc}
  \left\{  
    \begin{array}{cccc}
      P_\omega &&& \dim(P_\omega)=0\\
      \min(P_\omega) &&& 0 < \dim(P_\omega) < \dim(P) \\
      \texttt{EEP} &\quad&& P_\omega =P\\
    \end{array}
    \right.
\end{array}$$
where $\min(P_\omega)$ is the coordinate-wise minimum of all points in $P_\omega$ and \texttt{EEP} abbreviates \texttt{Exposes Entire Polytope}.  We remark that when a numerical oracle returns a vertex $v=\mathcal{O}_P(\omega)$, it also reveals that $\{x \in \mathbb{R}^n | \langle x, \omega \rangle \leq \langle v, \omega \rangle \}$ is a halfspace containing $P$. This fact is useful in finding a $V$-representation from an oracle \cite{OToV}. 

Given a polynomial
$$f=\sum_{\alpha \in \mathcal A} c_\alpha x_1^{\alpha_1}\cdots x_n^{\alpha_n}  \in \mathbb{C}[x_1,\ldots,x_n]\hspace{0.5 in} c_\alpha \neq 0, \mathcal A \subseteq \mathbb{N}^n, |\mathcal A|<\infty$$ its \emph{Newton polytope} $\textrm{New}(f)$ is the convex hull of $\mathcal A$. Motivated by language for polynomials, we say that $P$ is \emph{homogeneous} whenever $\mathcal O_P(1,1,\ldots,1) = \texttt{EEP}$ and define 
${\deg(P)}:=h_P(1,1,\ldots,1)$. The \emph{homogenization} of $P$ denoted {$\tilde P$} is the convex hull of $\{(x,\deg(P)-|x|) \big|x \in P\}$ where $|x|:=\sum_{i=1}^n x_i$. 

\begin{figure}[!htpb]
\begin{center}
\includegraphics[scale=0.2]{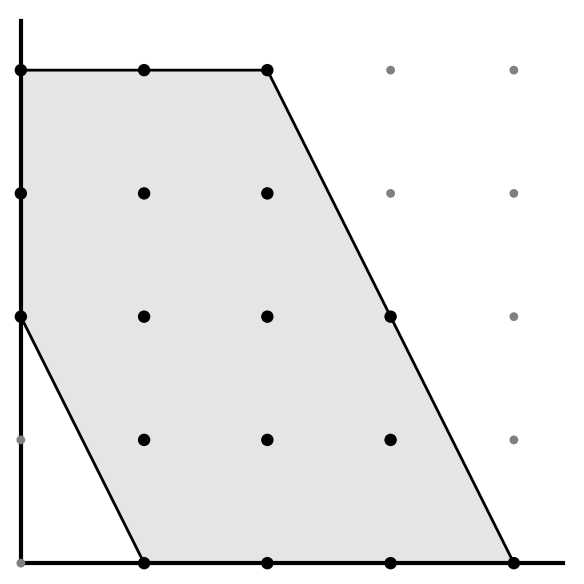} \hspace{0.5 in}
\includegraphics[scale=0.2]{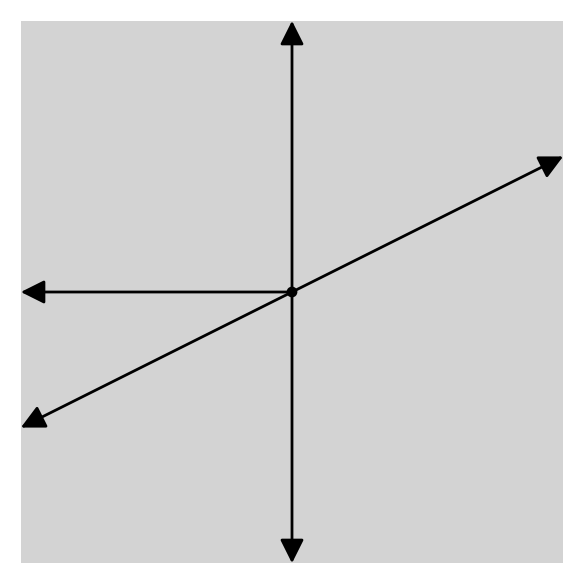}
\end{center}
\caption{A polytope and its corresponding normal fan}
\label{fig:polytope}
\end{figure}

 The \emph{(outer) normal fan} of $P$ is the polyhedral fan $\mathcal N(P)$ with cones $$C[\omega]=\{ \omega' \in \mathbb{R}^n | P_{\omega'} = P_\omega\}.$$
Figure \ref{fig:polytope} displays a polytope with vertices $(1,0),(4,0),(2,4),(0,4)$, and $(2,0)$ and its normal fan.
The fan has one zero-dimensional cone, five one-dimensional cones, and five two-dimensional cones.

\subsection{Some Tropical Geometry}
Newton polytopes are intimately related to tropical geometry. We only begin to touch on the topic here and encourage the reader to reference  \cite{MacSturm} for a more extensive treatment.

The tropicalization of a variety depends on the choice of a valuation on the base field involved (in our case $\mathbb{C}$). Relevant to our computations is the trivial valuation: $\nu(c)=0$ for all $c \in \mathbb{C}^*$. With this valuation, the tropicalization of a polynomial $$f=\sum_{\alpha \in \mathcal A} c_\alpha x_1^{\alpha_1}\cdots x_n^{\alpha_n}  \in \mathbb{C}[x_1,\ldots,x_n]\hspace{0.5 in} c_\alpha \neq 0, \mathcal A \subseteq \mathbb{N}^n, |\mathcal A|<\infty$$ 
is the map \begin{align*}\textrm{trop}(f):\mathbb{R}^n &\to \mathbb{R} \\
(\omega_1,\ldots,\omega_n) &\mapsto \max(\omega_i\cdot\alpha_i)
\end{align*} and the tropicalization of the hypersurface $V(f)$ is 
$$\textrm{trop}(V(f))=\{\omega \in \mathbb{R}^n| \text{ the maximum in }\textrm{trop}(f)(\omega)\text{ is attained at least twice}\}.$$
The tropicalization $V(f)$ is the codimension one fan of the normal fan of the Newton polytope of $f$. Moreover, $\textrm{trop}(V(f))$ is the locus of directions in which a numerical oracle does not return a vertex of $\textrm{New}(f)$. So, for example, if the polytope in Figure \ref{fig:polytope} is the Newton polytope of some hypersurface $V(f)$, then $\textrm{trop}(V(f))$ consists of the one-dimensional cones of the fan in Figure \ref{fig:polytope} along with origin.

The tropicalization of $V(I)$ for some ideal $I \subseteq \mathbb{C}[x_1,\ldots,x_n]$ is the intersection $$\textrm{trop}(V(I))=\bigcap_{f \in I} \textrm{trop}(V(f)).$$ Section \ref{Sec:Algorithms}  contains an algorithm to compute $\textrm{trop}(V(I))$ from projections.

\subsection{Numerical Algebraic Geometry}
Let $X\subseteq \mathbb{C}^N$ be an algebraic variety of dimension $k$ and degree $d$ appearing as an irreducible component of the zero set of a collection of polynomials $\mathcal F_X \subseteq \mathbb{C}[x_1,\ldots,x_n]$. For a generic $(N-k)$-dimensional linear space $\mathcal L \subseteq \mathbb{C}^N$, the intersection $X \cap \mathcal L$ is zero-dimensional and consists of $d$ points which are represented on a computer by a set $S$ of numerical approximations. The triple $(\mathcal F_X,\mathcal L, S)$ is called a \emph{witness set} for $X$ and is the fundamental data type in numerical algebraic geometry. The standard numerical method of \emph{homotopy continuation} quickly computes any witness set $(\mathcal F_X,\mathcal L', S')$ for $X$ from a precomputed witness set $(\mathcal F_X,\mathcal L,S)$ by numerically tracking the solutions $S$ along a homotopy from $S$ to $S'$ \cite{BertiniBook}. 
\begin{figure}[!htpb]
\centering
\includegraphics[scale=0.27]{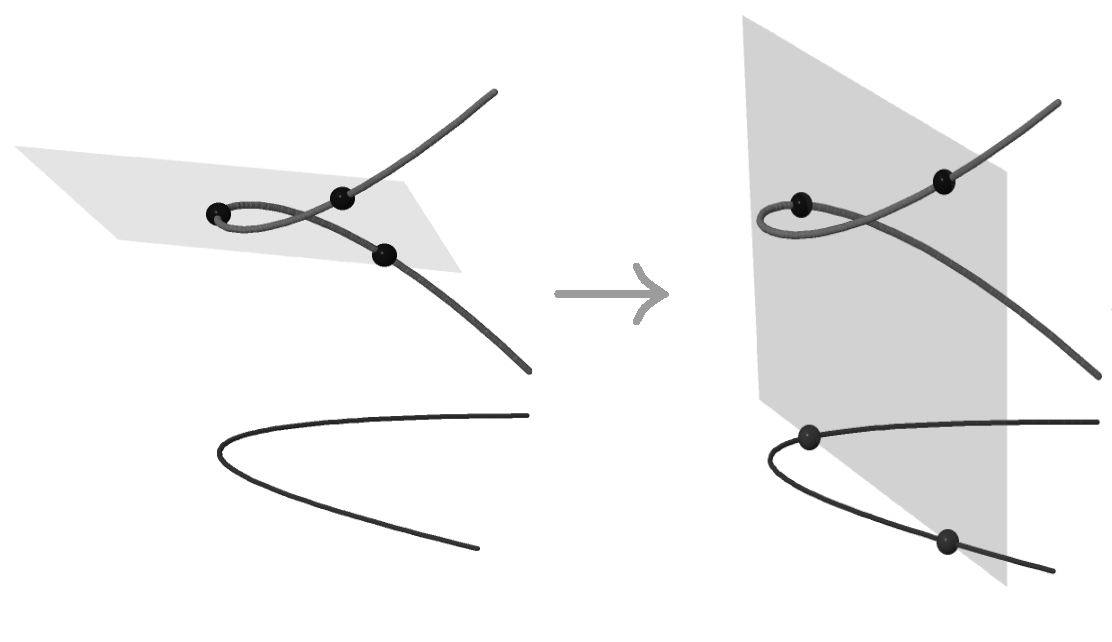}
\caption{Computing a witness set for a projection}
\label{fig:witnessProj}
\end{figure}
A major feature of numerical algebraic geometry is that we can compute witness sets for varieties without access to their equations. Let $X \subseteq \mathbb{C}^N$ be an irreducible and reduced component of a variety, $\pi:X \to \mathbb{C}^n$ a finite projection, and  $Y:=\overline{\pi(X)}$ the Zariski closure of its image. A witness set for $Y$ is encoded as a quadruple $(\mathcal F_X,\pi,\mathcal L',S')$ where $S'=\mathcal{L}' \cap Y$. Given a witness set $(\mathcal F_X,\mathcal L,S)$, we produce a witness set for $Y$ by performing a linear homotopy from the points in $S$ to the points $\pi^{-1}(S')=\pi^{-1}(\mathcal{L'}) \cap X$. For example, Figure \ref{fig:witnessProj} shows a witness set for the twisted cubic being
along with a witness set for the parabola coming from a projection of the twisted cubic. In the case that $\pi$ is not finite-to-one, we may still compute a witness set for the projection by replacing $X$ with $\hat X=X \cap L$ where $L$ is a generic linear space of dimension $\dim(X)-\dim(Y)$ so that $\dim(\hat X)=\dim(Y)$ and $\overline{\pi(\hat X)}=Y$. The fact that we can effectively compute witness sets for projections allows us to manipulate varieties which are presented as the images of maps since these are merely projections of graphs. Details on computing witness sets for projections can be found in \cite{WSP}.

\section{Algorithms}
\label{Sec:Algorithms}
\subsection{The HS-Algorithm}
Let $\mathcal H \subseteq \mathbb{C}^n$ be a degree $d$ hypersurface defined by $$f=\sum_{\alpha \in \mathcal A} c_\alpha x_1^{\alpha_1}\cdots x_n^{\alpha_n}  \in \mathbb{C}[x_1,\ldots,x_n]\hspace{0.5 in} c_\alpha \neq 0, \mathcal A \subseteq \mathbb{N}^n, |\mathcal A|<\infty$$ so that $\textrm{New}(f)$ is the convex hull of the points in $\mathcal A$. Let $\omega \in \mathbb{R}^n$ be a direction, $a,b \in (\mathbb{C}^*)^n$, and consider the family of lines $\mathcal L_t$ parametrized by
$$s \xrightarrow{{\bf L}_t} (\ell_{(t,1)}(s),\ldots,\ell_{(t,n)}(s))$$
where $\ell_{(t,i)}(s) = t^{\omega_i}(a_is-b_i)$. For any fixed $t$ value, $f({\bf L}_t)$ is a univariate polynomial in $s$ whose solutions $p(t):=\{p_{1}(t),\ldots,p_{d}(t)\}$ in $\mathbb{C}_s$ correspond to intersection points of $\mathcal H$ and $\mathcal L_t$. We may write $f({\bf L}_t)$ as  
\begin{align*}
f(\ell_{(t,1)},\ldots,\ell_{(t,n)})&=\sum_{\alpha \in \mathcal A} c_\alpha [t^{\omega_1}(a_1s-b_1)]^{\alpha_1}\cdots [t^{\omega_n}(a_ns-b_n)]^{\alpha_n}\\
&=\sum_{\alpha\in \mathcal A}t^{\langle \omega, \alpha\rangle}(a_1s-b_1)^{\alpha_1}\cdots (a_ns-b_n)^{\alpha_n}.
\end{align*}
As $t \to \infty$, the terms $\mathcal A_\omega$ corresponding to points of $\mathcal A$ which maximize $\langle \omega, \alpha \rangle$ will dominate the behavior of the zeros and so the solutions $p(t)$ will converge to those of 
$$f_\omega(\ell_{(t,1)},\ldots,\ell_{(t,n)}):= \sum_{\alpha \in \mathcal A_{\omega}} c_\alpha(a_1s-b_1)^{\alpha_1}\cdots (a_ns-b_n)^{\alpha_n}.$$
If $\mathcal A_{\omega}=\{\beta\}$ then $f_\omega$ is a monomial and so $f_\omega({\bf L}_t)$ has roots $\gamma_i:=b_i/a_i$ where $\gamma_i$ occurs with multiplicity $\beta_i$. 
If $|\beta|:=\sum \beta_i$ is less than $d$, then there are $\beta_\infty:=d-|\beta|$ points which have diverged towards infinity. One can see this by observing that if we began with the homogenization $F$ of $f$, this would be the exponent of the homogenizing variable in the term $F_\omega$. 

If $\omega$ exposes the entire polytope defined by $\mathcal A$, then the roots $p(t)$ remain constant as $t$ varies since all $f({\bf L}_t)$ are all scalar multiples of each other. 
 
If $\omega$ exposes a proper non-trivial subset of $\mathcal A$, then there is more than one term in $f_\omega$, or equivalently $\omega \in \textrm{trop}(V(f))$. The terms of $f_\omega$ will have a common factor of $\prod_{i=1}^n (a_is-b_i)^{m_i}$ where the vector $m$ is the coordinate-wise minimum of the points in $\mathcal A_\omega$. Therefore,  $m_i$ roots will converge to $\gamma_i$ and $m_\infty:=\min_{\beta \in \mathcal A_\omega} \left(d-|\beta|\right)$ points will diverge to infinity. All other roots will converge somewhere else in $\mathbb{C}$. 

These observations give rise to the HS-algorithm.\\ \\ 
$\begin{array}{l}
  \left[ 
    \begin{array}{l}
        \\
         \\
        \\
         \\
        \\
         \\
        \\
         \\
        \\
         \\
        \\
         \\
        \\
         \\
        \\
         \\
        \\
        \\
        \\
    \end{array}
    \right.
\end{array}$
\begin{minipage}{\textwidth}{}
\begin{alg}{\texttt{HS-Algorithm}}
\label{HSAlgorithm}\\
{\bf Input:}
\begin{itemize}
\item A witness set $W$ for a hypersurface $\mathcal H \subseteq \mathbb{C}^n$
\item A direction $\omega \in \mathbb{R}^n$
\end{itemize}
{\bf Output:}
\begin{itemize}
\item $\mathcal O_{\widetilde{\textrm{New}(\mathcal H)}}(\omega)$
\end{itemize}
{\bf Steps:}
\begin{enumerate}
\item Pick random $a,b \in \mathbb{C}^n$ and construct $\{\ell_{(t,i)}\}_{i=1}^n$ described above
\item Track the witness points in $W$ to the intersection $\mathcal H \cap \mathcal L_1$
\item Initialize vector $\beta={\bf 0} \in \mathbb{N}^{n+1}$
\item Track the witness points in $\mathcal H \cap \mathcal L_t$ from $t=1$ toward $\infty$
\item If none of the solutions move, return \texttt{EEP}
\item If a solution has converged, stop tracking it
\begin{itemize}
\item If it has converged to some $\gamma_i$ increment $\beta_i$ by one
\end{itemize}
\item If a solution has diverged  increment $\beta_\infty$ by one
\item If all solutions have converged or diverged, return $\beta=(\beta_1,\ldots,\beta_n,\beta_\infty)$
\end{enumerate}
\end{alg}
\end{minipage}

\subsection{Tropical Membership}
\label{Section:Tropical}

Motivated by the results of Bieri and Groves in \cite{BieriGroves}, Hept and Theobald in \cite{HeptTheobald} investigated how to write $\textrm{trop}(V(I))$ as an intersection of tropical hypersurfaces coming from projections. The following is a consequence of the proof of Theorem 1.1 in \cite{HeptTheobald}.
\begin{theorem}
 If $I\subseteq \mathbb{C}[x_1,\ldots,x_n]$ is an $m$-dimensional prime ideal, and $\{\pi_i\}_{i=0}^{n-m}$ are generic projections, 
 $$\textrm{trop}(V(I)) =  \bigcap_{i=0}^{n-m} \pi_i^{-1}(\pi_i(\textrm{trop}(V(I)))$$ 
 where each $\pi_i^{-1}(\pi_i(\textrm{trop}(V(I))))$ is a tropical hypersurface. 
\label{TropicalProjection}
\end{theorem}

Unfortunately, coordinate projections are not always generic and it is possible that we have a proper containment $$\bigcap_{J \subseteq [n]: \text{codim}(\pi_J(V(I)))=1} \pi_J^{-1}(\pi_J(V(I))) \subsetneq \textrm{trop}(V(I)) $$  where $\pi_J$ denotes projection onto the coordinates $\{x_j\}_{j \in J}$ . When this is the case, it is necessary to take more general projections of $\textrm{trop}(V(I))$ as in Example \ref{ex:trop}.
 
The following algorithm is a consequence of Theorem \ref{TropicalProjection} and the HS-algorithm.
\\
$\begin{array}{l}
  \left[ 
    \begin{array}{l}
        \\
        \\
         \\
        \\
        \\
         \\
        \\
         \\
        \\
         \\
        \\
         \\
        \\
        \\
        \\
        \\
        \\
        \\
        \\
        \\
        \\
    \end{array}
    \right.
\end{array}$
\begin{minipage}{0.9\textwidth}{}
\begin{alg}{Tropical Membership}
\label{TropicalMembership}\\
{\bf Input:}
\begin{itemize}
\item An $m$-dimensional variety $X=V(I)\subseteq \mathbb{C}^n$
\item A direction $\omega \in \mathbb{R}^n$
\end{itemize}
{\bf Output:}
\begin{itemize}
\item \texttt{true} if $\omega \in \textrm{trop}(X)$ and \texttt{false} otherwise.
\end{itemize}
{\bf Steps:}
\begin{enumerate}
\item First, replace $I$ with its image under a random monomial map $\varphi(x_i)=\prod_{j=1}^n(x_j)^{A_{(i,j)}}$ so that the coordinate projections of $V(I)$ are generic. Simultaneously replace $\omega$ with $A^{-1}\omega$.
\item Compute a witness set $W$ for $X$.
\item For all coordinate projections $\{\pi_J\}_{J \subseteq [n]}$ with $|J|=n-m-1$, use $W$ to compute a witness set $W_J$ for $\pi_J(X)$.
\item Using each witness set, run the HS-algorithm on $\pi_J(X)$ in the direction $\pi(\omega)$. 
\item If, for each such projection, the HS-algorithm observes convergence of all solutions, but $|\beta| \neq \deg(\pi(X))$, return $\texttt{true}$, otherwise return $\texttt{false}$.
\end{enumerate}
\end{alg}
\end{minipage}

We remark that the monomial change of coordinates involved in step $(1)$ may enlarge the degree of $X$ thus making the computation of a witness set more difficult. It is, however, often the case that the coordinate projections of $\textrm{trop}(X)$ are already general without any monomial change of coordinates. Moreover, when this is not the case, the algorithm can only yield false positives.

Theorem 8 of \cite{NP}, gives an analysis of the convergence of the HS-algorithm whenever $\omega$ exposes a vertex. We provide an analogous result for the case where $\omega \in \textrm{trop}(V(f))$ and thus $f_\omega$ correpsonds to a positive dimensional face of $\textrm{New}(f)$  with monomial support $\mathcal F\subseteq \mathcal A$. 

We first introduce notation. Let $x^m$ be the common monomial factor of $f_\omega$ when $\omega$ exposes a positive dimensional face and write $f_\omega = x^m\cdot g_\omega$. Also define $k_{\omega,a,b}$ to be the constant appearing in
$$g_\omega({\bf L}_t) =t^{h(\omega)}\cdot  k_{\omega,a,b} \cdot \prod_{i=1}^{d-|m|-m_\infty} ( s-\tau_i)^{\beta_{\tau_i}}, \hspace {0.2 in}k_{\omega,a,b} \in \mathbb{C}.$$
Set $a_{\min}:=\min\{1,|a_i|:i=1,\ldots,n\}$, $a_{\max}:=\max\{1,|a_i|: i=1,\ldots,n\}$. Define,
$$C:=\frac{\max(\{|c_\alpha|:\alpha \in \mathcal A\})}{k_{\omega,a,b}}$$
$$\gamma_\tau:=\min\left\{a_{\min},\frac{1}{2}|\tau-\rho|\text{ such that }\rho \in \left\{\frac{b_i}{a_i}\right\}_{i=1}^n \cup \{\tau_i\}_{\tau_i \neq \tau} \right\}$$
$$\Gamma_\tau:=\max\left\{\frac{2}{a_{\max}},|\tau-\rho|\text{ such that }\rho \in \left\{\frac{b_i}{a_i} \right\}_{i=1}^n \right\}$$
$$d_\omega = \min \{h(\omega)-\langle \omega , \alpha\rangle :\alpha \in \mathcal F^c\}.$$

Since it is enough to observe convergence of some path in the HS-algorithm to a point in $\mathbb{C}$ other than $\frac{b_i}{a_i}$ for some $i=1,\ldots,n$, we analyze the convergence rate for such paths only.
\begin{theorem}
Suppose $\omega \in \textrm{trop}(V(f))$.
Let $p(t)$ be a path of the HS-algorithm converging to $\tau\not\in \left\{\frac{b_i}{a_i}\right\}_{i=1}^n$ as $t \to \infty$ and let $\beta$ be the number of such paths converging to $\tau$. Let $t_1 \geq 0$ be a number such that if $t >t_1$ then 
$|p(t) - \tau| \leq \gamma_\tau$.
Then for all $t > t_1$ 
$$|p(t)-\tau|^\beta \leq t^{-d_\omega} \cdot C \cdot |\mathcal F^c| \cdot \left(\frac{a_{\max}}{a_{\min}} \left(1+\frac{\Gamma_\tau}{\gamma_\tau}\right)\right)^d.$$
\end{theorem}
{\it Proof of Theorem 2:}\\
Writing
$$f({\bf L}_t)=f_\omega({\bf L}_t)+\sum_{\alpha \in \mathcal F^c} c_\alpha {\bf L}_t^\alpha$$ we have
$$f({\bf L}_t) = t^{h(\omega)}f_\omega(as-b)+\sum_{\alpha \in \mathcal F^c} c_\alpha (as-b)^{\alpha}t^{\langle \alpha, \omega \rangle}$$ but since $f(p(t))\equiv 0$ we have
\begin{align*}
0=t^{h(\omega)}f_\omega(ap(t)-b)&+t^{h(\omega)}\sum_{\alpha \in \mathcal F^c} c_{\alpha}(ap(t)-b)^\alpha t^{\langle \alpha, \omega \rangle -h(\omega)}\\
\implies |f_\omega(ap(t)-b)| &= \bigg|\sum_{\alpha \in \mathcal F^c} c_\alpha(ap(t)-b)^\alpha  t^{\langle \alpha, \omega \rangle -h(\omega)}\bigg |\\
 &\leq t^{-d_\omega} \sum_{\alpha \in \mathcal F^c} |c_\alpha|\cdot|(ap(t)-b)^\alpha |.
 \end{align*}
However, since $f_\omega(as-b)=(as-b)^m k_{\omega,a,b}\prod (s-\tau_i)^{\beta_i}$ we divide through by $k_{\omega,a,b}$ so that 
$$|(as-b)^m|\cdot \prod |(s-\tau_i)|^{\beta_i} \leq  t^{-d_\omega}\cdot C \cdot  \sum_{\alpha \in \mathcal F^c} |(ap(t)-b)^\alpha |.$$
We now bound the right-hand-side summands
\begin{align*}
|a_jp(t)-b_j| &= |a_j|\cdot |p(t)-b_j/a_j|\\
& \leq a_{\max}|p(t)-\tau+\tau-b_j/a_j|\\
&\leq a_{\max}(\gamma_\tau+\Gamma_\tau)
\end{align*} so that 
$$|(ap(t)-b)^\alpha| \leq |a_{\max}(\gamma_\tau+\Gamma_\tau)|^d$$
 because $2 \leq a_{\max} \Gamma_\tau$ and $|\alpha| \leq d$.
Now substitution yields 
\begin{align*}
|(ap(t)-b)^m|\cdot \prod |(p(t)-\tau_i)|^{\beta_i} 
&\leq t^{-d_\omega}\cdot C \cdot  \sum_{\alpha \in \mathcal F^c}  |a_{\max}(\gamma_\tau+\Gamma_\tau)|^d\\
&\leq t^{-d_\omega}\cdot C \cdot  |\mathcal F^c|\cdot   |a_{\max}(\gamma_\tau+\Gamma_\tau)|^d .
\end{align*}
Let $g(t):=|(ap(t)-b)^m|\cdot \prod |(p(t)-\tau_i)|^{\beta_i}$. We now bound the size of the factors of $g(t)$ other than $(p(t)-\tau)^{\beta}$.
\begin{align*}
|p(t)a_j-b_j| = |a_j|\cdot |p(t)-b_j/a_j|&=|a_j| \cdot |p(t)-\tau+\tau-b_j/a_j|\\
&\geq a_{\min} \cdot \big| |\tau- b_j/a_j|-|p(t)-\tau|\big|\\
&\geq a_{\min}|2\gamma_\tau-\gamma_\tau| = a_{\min} \gamma_\tau
\end{align*} and similarly $|p(t)-\tau_j| \geq \gamma_\tau \geq a_{\min}\gamma_\tau$ for $ \tau_j \neq \tau$. Since $a_{\min} \gamma_\tau \leq 1$ and $\beta \leq d$
$$\bigg|\frac{g(t)}{(p(t)-\tau)^\beta}\bigg| \geq (a_{\min}\gamma_\tau)^{d-\beta}\geq (a_{\min}\gamma_\tau)^{d}$$ so
\begin{align*}
|p(t)-\tau|^{\beta} &\leq t^{-d_\omega} \cdot C \cdot |\mathcal F^c| \cdot |a_{\max}(\gamma_\tau+\Gamma_\tau)|^d\cdot \frac{1}{(a_{\min}\gamma_\tau)^{d}}\\
&\leq t^{-d_\omega} \cdot C \cdot |\mathcal F^C| \cdot \left( \frac{a_{\max}}{a_{\min}} \left(1+ \frac{\Gamma_\tau}{\gamma_\tau}\right)\right)^d. \hspace{1 in}  \square
\end{align*}

\section{Functionality}
\label{Section:Functionality}
There are four main user functions in {\bf NumericalNP.m2}. The first three implement the HS-algorithm and the last implements the tropical membership algorithm.

Function \ref{witnessForProjection}, computes a witness set for the image of an irreducible and reduced variety $X \subseteq \mathbb{C}^N$  under a projection $\pi: \mathbb{C}^N \to \mathbb{C}^n$. \\ \\
$\begin{array}{c}
  \left[ 
    \begin{array}{c}
        \\
         \\
        \\
         \\
        \\
         \\
        \\
         \\
        \\
        \\
         \\
        \\
    \end{array}
    \right.
\end{array}$
\begin{minipage}{0.90\textwidth}{}
\begin{funct}{\texttt{witnessForProjection}}
\label{witnessForProjection}\\
{\bf Input:}
\begin{itemize}
\item I: Ideal defining $X \subseteq \mathbb{C}^N$
\item ProjCoord: List of coordinates which are forgotten by $\pi$
\item OracleLocation (option): Path in which to create witness files
\end{itemize}
{\bf Output:} A subdirectory \texttt{/OracleLocation/WitnessSet} containing
\begin{itemize}
\item witnessPointsForProj: Preimages of witness points of $\overline{\pi(X)}$
\item projectionFile: List of coordinates in ProjCoord
\item equations: List of equations defining $X' \subseteq X$ such that $\pi|_{X'}$ is finite and $\overline{\pi(X')}=\overline{\pi(X)}$
\end{itemize}
\end{funct}
\end{minipage}

Function \ref{witnessToOracle}, \texttt{witnessToOracle},  creates all necessary {\bf Bertini} files to track the witness set $\mathcal H \cap \mathcal L_t$ as $t \to \infty$ for any $\omega \in \mathbb{R}^n$. These files treat $\omega$ as a parameter so that the user only needs to produce these files once. \\ \\
$\begin{array}{c}
  \left[ 
    \begin{array}{c}
        \\
         \\
        \\
         \\
        \\
         \\
        \\
        \\
        
        \\
        \\
        \\
        \\
        \\
    \end{array}
    \right.
\end{array}$
\begin{minipage}{0.95\textwidth}{}
\begin{funct}{\texttt{witnessToOracle}}
\label{witnessToOracle}\\
{\bf Input:}
\begin{itemize}
\item  OracleLocation: Path containing the directory \texttt{/WitnessSet}
\end{itemize}
{\bf Optional Input:}
\begin{itemize}
\item \texttt{PointChoice}: Prescribes $a$ and $b$ explicitly (see Algorithm \ref{HSAlgorithm})
\item \texttt{TargetChoice}: Prescribes targets $b_i/a_i$
\item \texttt{NPConfigs}: List of {\bf Bertini} path tracking configurations
\end{itemize}
{\bf Output:} 
\begin{itemize}
\item A subdirectory \texttt{/OracleLocation/Oracle} containing all necessary files to run the homotopy described in Algorithm \ref{HSAlgorithm}.
\end{itemize}
\end{funct}
\end{minipage}

Function \ref{witnessToOracle} by default chooses  $a,b \in \mathbb{C}^n$ such that $\gamma_i :=a_i/b_i$ are $n$-th roots of unity. One may choose to either specify $a$ and $b$ (\texttt{PointChoice}), or $\gamma_i:=a_i/b_i$ (\texttt{TargetChoice}) or request that these choices are random. When random, the function ensures that the points $\gamma_i$ are far from each other so that convergence to $\gamma_i$ is easily distinguished from convergence to $\gamma_j$.  
{\bf Bertini} is called to track the solutions in \texttt{/OracleLocation/WitnessSet} to points $\overline{\pi(X)} \cap \mathcal L_1$. These become start solutions to the homotopy described in Algorithm \ref{HSAlgorithm} with parameters $\omega$ and $t$. There are many numerical choices for {\bf Bertini}'s native pathtracking algorithms which can be specified via \texttt{NPConfigs}.\\ \\ 
$\begin{array}{c}
  \left[ 
    \begin{array}{c}
        \\
         \\
        \\
         \\
        
         \\
        \\
        \\
        \\
        \\
        \\
        \\
        \\
        
        \\
        \\
    \end{array}
    \right.
\end{array}$
\begin{minipage}{0.94\textwidth}{}
\begin{funct}{\texttt{oracleQuery}}
\label{oracleQuery}

{\bf Input:}

\begin{itemize}
\item   OracleLocation (Option): Location containing the directory \texttt{/Oracle}
\item  $\omega$: A vector in $\mathbb{R}^n$
\end{itemize}

{\bf Optional Input:}\\
$-$ \texttt{Certainty} \hspace{0.05 in}  $-$ \texttt{Epsilon} \hspace{0.05 in}  $-$ 
  \texttt{MinTracks} \hspace{0.05 in}  $-$ 
 \texttt{MaxTracks} \hspace{0.05 in}
 
   $-$ 
 \texttt{StepResolution} \hspace{0.05 in}  $-$ 
 \texttt{MakeSageFile}

{\bf Output:}

\begin{itemize}
\item $\mathcal O_{\widetilde{\textrm{New}(\mathcal H)}}(\omega)$ or \texttt{Reached MaxTracks}
\item  A subdirectory \texttt{/OracleLocation/OracleCalls/Call\#} containing
\begin{itemize}
\item[--]  \texttt{SageFile}: Sage code animating the paths $p(t)$
\item[--]  \texttt{OracleCallSummary}: a human-readable file summarizing the results
\end{itemize}
\end{itemize}
\end{funct}
\end{minipage}

The fundamental function, \texttt{oracleQuery}, runs the homotopy in the HS-algorithm, monitors convergence, and outputs the result of the numerical oracle.

To monitor convergence of solutions $p(t)$ we track $t \to \infty$ in discrete steps. The option \texttt{StepResolution} specifies these $t$-step sizes. In each step, for each path $p_{i}(t)$, a numerical derivative is computed to determine convergence or divergence of the solution. If the solution is large and the numerical derivative exceeds $10^{\texttt{Certainty}}$ in two consecutive steps the path is declared to diverge, and if the numerical derivative is below $10^{-\texttt{Certainty}}$ in two consecutive steps the point is declared to converge. If a converged point is at most \texttt{Epsilon} from some $\gamma_i$, then the software deems that it has converged to $\gamma_i$. When a point is declared to converge or diverge, it is not tracked further.
\begin{figure}[!htpb]
\centering
\includegraphics[scale=0.18]{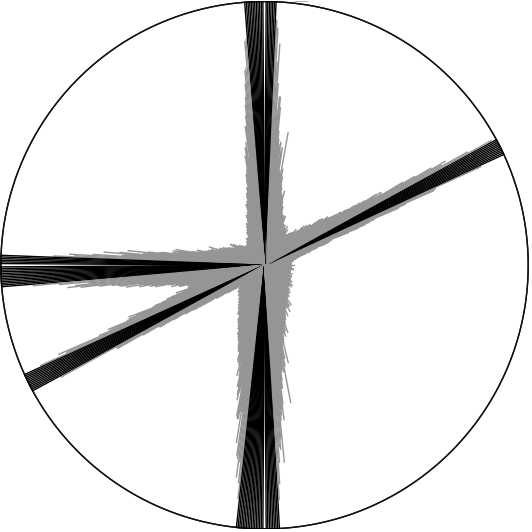} \hspace{0.2 in}
\includegraphics[scale=0.18]{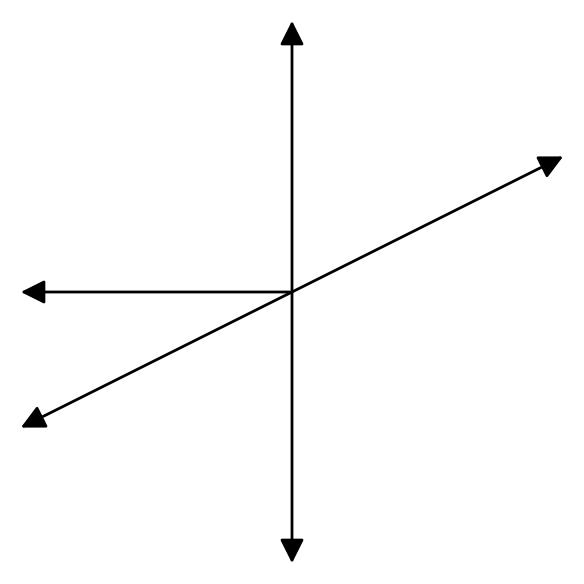}
\hspace{0.2 in}
\includegraphics[scale=0.18]{polytope.png}
\caption{ Convergence rates (Left) of different $\omega$ for queryOracle on hypersurface with Newton polytope (Right) and normal fan (center)}
\label{koosh}
\end{figure}
The option \texttt{MaxTracks} allows the user to specify how long to wait for convergence of the paths $p(t)$.  Figure \ref{koosh} shows the Newton polytope of a plane sextic (see Example \ref{ex:example}) as well as the convergence rate of the algorithm on different directions $\omega \in S^1$: the length of each grey ray is proportional to the number of steps required to observe convergence and the black rays indicate that this convergence was not observed within the limit specified by \texttt{MaxTracks}. We include the image of the tropicalization of this curve to illustrate how the convergence rate involved in the HS-algorithm slows as $\omega$ approaches directions in the tropical variety. Nonetheless, when $\omega$ is {\it precisely} in the tropical variety, the runtime is actually quite small, evident in the small gaps in the tropical directions of Figure \ref{koosh}. 

One may also specify \texttt{MinTracks} which indicates the step at which convergence begins to be monitored. The option to create a Sage \cite{sage} animation (see Figure \ref{SAGE}) of the solution paths helps the user recognize pathological behavior in the numerical computations and fine-tune parameters such as \texttt{Certainty}, \texttt{StepResolution}, or \texttt{Epsilon} accordingly.
\begin{example}
\label{ex:example}
Consider the curve in $X \subseteq \mathbb{C}^3$ defined by $$I=\langle xyt-(x-y-t)^2+3x+t,x+y^2+t^2 \rangle\subseteq \mathbb{C}[x,y,t]$$ and let $\pi$ be the projection forgetting the $t$ coordinate. The following {\bf Macaulay2} code computes a witness set for $\mathcal C:=\overline{\pi(X)}$, prepares oracle files for the HS-algorithm and then runs the HS-algorithm in the direction $(3,2)$. The software returns $\{2,4,0\}$ indicating that $\textrm{New}(\overline{\pi(X)})_{(3,2)}=(2,4)$.
{\small
\begin{verbatim}
i1: loadPackage("NumericalNP");
i2: R=CC[x,y,t];
i3: I=ideal(x*y*t-(x-y-t)^2+3*x+t,x+y^2+t^2);
i4: witnessForProjection(I,{2},OracleLocation=>"Example");
i5: witnessToOracle("Example") ;
i6: time oracleQuery({3,2},OracleLocation=>"Example",MakeSageFile=>true)
     -- used 0.178448 seconds
o6: {2,4,0} 
\end{verbatim}}
The full Newton polytope of $\overline{\pi(X)}$ is displayed in Figure \ref{koosh} and snapshots of the Sage animation created by \texttt{queryOracle} are shown in Figure \ref{SAGE}. There, the circles are centered at $\gamma_1=1$ and $\gamma_2=-1$ and have radius \texttt{epsilon}.
\begin{figure}[!htpb]
\centering
\includegraphics[scale=0.26]{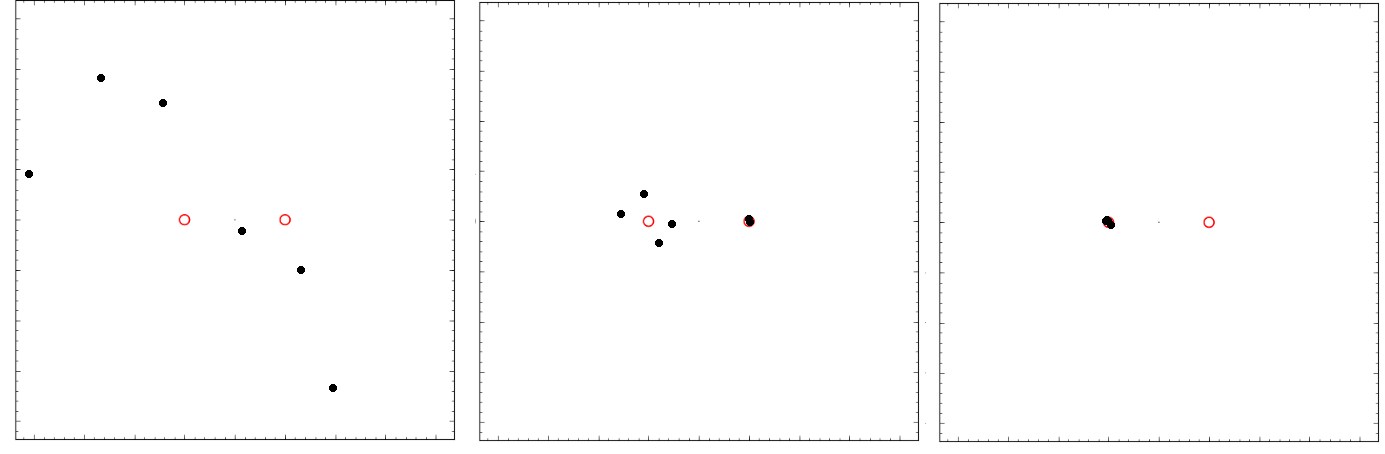}
\caption{Three snapshots of Sage animation from example with viewing window $[-4,4]^2$}
\label{SAGE}
\end{figure}
\end{example}
$\begin{array}{c}
  \left[ 
    \begin{array}{c}
        \\
         \\
        \\
         \\
        \\
         \\
        \\
        \\
        \\
        \\
        \\
        \\
        \\
        \\
    \end{array}
    \right.
\end{array}$
\begin{minipage}{0.94\textwidth}{}
\begin{funct}{\texttt{tropicalMembership}}
\label{tropicalMembership}

{\bf Input:}

\begin{itemize}
\item   $I:$Ideal defining $X \subseteq \mathbb{C}^n$
\item  $\omega$: A vector in $\mathbb{R}^n$
\end{itemize}

{\bf Optional Input:}\\
$-$ \texttt{Certainty} \hspace{0.05 in}  $-$ \texttt{Epsilon}\hspace{0.05 in}  \hspace{0.05 in}  $-$ 
  \texttt{MinTracks} \hspace{0.05 in}  $-$ 
 \texttt{MaxTracks} \hspace{0.05 in}
 
   $-$ 
 \texttt{StepResolution} \hspace{0.05 in}  $-$ 
 \texttt{MakeSageFile}

{\bf Output:}

\begin{itemize}
\item A list of oracle queries of $\pi(X)$ in the direction $\pi(\omega)$ where $\pi$ runs through all coordinate projections such that $\pi(X)$ is a hypersurface.
\item \texttt{true} if all oracle queries exposed positive dimensional faces and \texttt{false} otherwise
\end{itemize}
\end{funct}
\end{minipage}

The fourth function \texttt{tropicalMembership} computes a witness set for each coordinate projection $\pi(V(I))$ whose image is a hypersurface. The algorithm subsequently checks that \texttt{oracleQuery} indicates that $\pi(\omega) \in \textrm{trop}(\pi(V(I))$. If this is true for each coordinate projection, the algorithm returns \texttt{true} and otherwise returns \texttt{false}. The options fed to \texttt{tropicalMembership} are passed along to \texttt{oracleQuery}.
\\ 
\begin{example}
\label{ex:trop}
Example 4.2.11 in \cite{Chan} gives two tropical space curves which are different, yet have the same tropicalized coordinate projections. We depict these in Figure \ref{fig:badtrop} and illustrate this behavior with our software.
\begin{figure}
\centering
\includegraphics[scale=0.3]{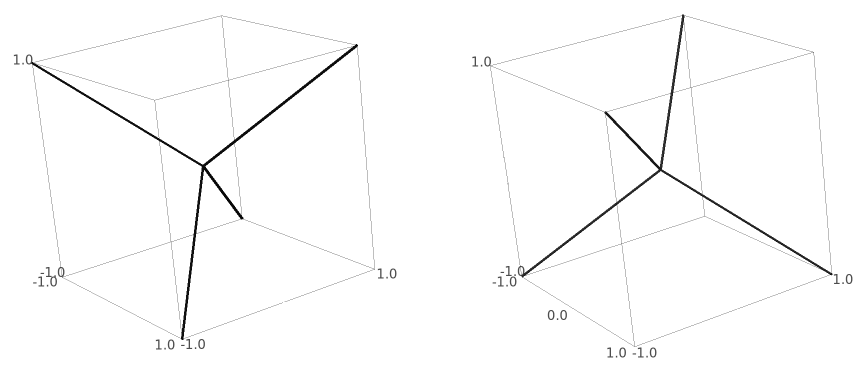} 
\caption{Two tropical space curves with the same tropical coordinate projections}
\label{fig:badtrop}
\end{figure}
{\small
\begin{verbatim}
i1 : loadPackage("NumericalNP");
i2 : R=QQ[x,y,z];
i3 : I=ideal {x*z+4*y*z-z^2+3*x-12*y+5*z,x*y-4*y^2+y*z+x+2*y-z};
i4 : J=ideal{x*y-3*x*z+3*y*z-1,3*x*z^2-12*y*z^2+x*z+4*y*z+5*z-1};
i5 : I==J
o5 = false
i6 : directions:={{1,1,1},{1,1,-1},{1,-1,1},
{1,-1,-1},{-1,1,1},{-1,1,-1},{-1,-1,1},{-1,-1,-1}};
i7 : apply(directions,d->tropicalMembership(J,d))
o7 = {true, true, true, true, true, true, true, true}
i8 : apply(directions,d->tropicalMembership(I,d))
o8 = {true, true, true, true, true, true, true, true}
\end{verbatim}}
\end{example}
Every projection of every vertex of the $\{\pm 1\}^3-$cube is in the tropicalization of the corresponding projection of $V(I)$ and $V(J)$. Nonetheless, the tropicalizations of $V(I)$ and $V(J)$ are disjoint subsets of the vertices of the cube. This exemplifies that an output of \texttt{true} from \texttt{tropicalMembership} is not a certification of membership of the tropical variety as we cannot {\it a priori} decide whether or not our coordinate projections are generic. When the projection is sufficiently generic, \texttt{tropicalMembership} will return \texttt{true} if and only if $\omega \in \textrm{trop}(V(I))$. \\
\\
{\it Example \ref{ex:trop} (continued):}
Consider monomial change of coordinates $\varphi$ given by $\varphi(x)=xyz,\varphi(y)=y$, and $\varphi(z)=z$. Let $I'$ be the extension of $I$ under this map so that $\varphi:\mathbb{C}[x,y,z]/I \to \mathbb{C}[x,y,z]/I'$ and $\varphi^*:V(I') \to V(I)$. By Corollary 3.2.13 of \cite{MacSturm} we have that$$\textrm{trop}(V(I)) = \textrm{trop}(\varphi)(\textrm{trop}(V(I'))), \hspace{0.2 in} \textrm{trop}(\varphi)={\small
\begin{bmatrix}
1 & 1 & 1 \\
0 & 1 & 0 \\
0 & 0 & 1
\end{bmatrix}}={\small \begin{bmatrix}
1 & -1 & -1 \\
0 & 1 & 0 \\
0 & 0 & 1
\end{bmatrix}^{-1}.}
$$
In other words, $\varphi$ induces a linear transformation $\textrm{trop}(\varphi)$ on tropical varieties. This transformation is generic in the sense of Theorem \ref{TropicalProjection} and \texttt{tropicalMembership} is able to distinguish $\textrm{trop}(V(I))$ from $\textrm{trop}(V(J))$.
{\small
\begin{verbatim}
i9 : I'=ideal apply(I_*,f->sub(f,{x=>x*y*z,y=>y,z=>z}));
i10 : J'=ideal apply(J_*,f->sub(f,{x=>x*y*z,y=>y,z=>z}));
i11 : directions'=apply(directions,d->{d#0-d#1-d#2,d#1,d#2})
o11 = {{-1, 1, 1}, {1, 1, -1}, {1, -1, 1}, {3, -1, -1},
      {-3, 1, 1}, {-1, 1, -1}, {-1, -1, 1}, {1, -1, -1}}
i12 : apply(directions',d->tropicalMembership(I',d))
o12 = {false, true, true, false, true, false, false, true}
i13 : apply(directions',d->tropicalMembership(J',d))
o13 : {true, false, false, true, false, true, true, false}
\end{verbatim}
}

\section{Applications}
\label{Section:Applications}

\subsection{The L\"uroth Polytope}
A \emph{L\"uroth quartic} is a plane quartic which interpolates the ten intersection points of a configuration of five lines in the plane.
The set of all L\"uroth quartics is a rational hypersurface $\mathcal H$ of degree $54$ in the $15$ coefficients $\{q_{(i,j,k)}\}_{i+j+k=4}$ of a plane quartic called the \emph{L\"uroth hypersurface} and is defined by a single homogeneous polynomial $f$. This hypersurface is an invariant of $GL(\mathbb{C}^3)$ and so the permutation subgroup $G=S_3$ acts on the vertices of $\textrm{New}(\mathcal H)$ by permuting the three indeterminants of a homogeneous quartic. A face of $\textrm{New}(\mathcal H)$ was found in \cite{NP}. Using our software, we have rediscovered that $\textrm{New}(\mathcal H)$ is $12$-dimensional and have, so far, found $1713$ vertices, belonging to $1,1,28,$ and $271$ orbits of sizes $1,2,3,$ and $6$ respectively. 

Querying the oracle in the coordinate directions bounds $\textrm{New}({\mathcal H})$ in a box. These bounds are given by $q_{(4,0,0)} \in [0,18]$, $q_{(3,1,0)} \in [0,24]$, $q_{(2,2,0)} \in [0,28]$, and $q_{(2,1,1)} \in [0,32]$ up to permutation of the indices. 

Up-to-date computations regarding the L\"uroth invariant as well as the package {\bf NumericalNP.m2} can be found at the authors webpage \cite{taylor}.

\subsection{Algebraic Vision Tensor}
\label{sec:algVision}
The multiview variety $X$ of a pinhole camera and a two slit camera is a hypersurface in the space of  $3 \times 2 \times 2$ tensors given by the image of twelve particular minors of 
$$\begin{bmatrix}
A&B&C
\end{bmatrix}
={\small
\begin{bmatrix}
a_{1,1} & a_{1,2} & a_{1,3} & b_{1,1} & b_{1,2} & c_{1,1} & c_{1,2} \\
a_{2,1} & a_{2,2} & a_{2,3} & b_{2,1} & b_{2,2} & c_{2,1} & c_{2,2} \\
a_{3,1} & a_{3,2} & a_{3,3} & b_{3,1} & b_{3,2} & c_{3,1} & c_{3,2} \\
a_{4,1} & a_{4,2} & a_{4,3} & b_{4,1} & b_{4,2} & c_{4,1} & c_{4,2} 
\end{bmatrix}.}
$$ We consider $X$ as a subvariety of $\mathbb{P}^{11}$ given by the image of 
$$F:\mathbb{C}^{28} \to \mathbb{C}^{12} \hspace{ 0.8 in} 
[A B C] \xrightarrow{F} \{f_{i,j,k}\}_{i \in \{1,2,3\},j,k \in \{1,2\}}$$
where $f_{i,j,k}$ is the minor not involving columns $a_i,b_j,$ and $c_k$. This map has $17$-dimensional fibres so \texttt{witnessForProjection} automatically slices $\mathbb{C}^{28}$ with $17$ hyperplanes to compute a witness set for $X$ which shows that $\deg(X)=6$. Therefore, its defining polynomial has an {\it a priori} upper bound of $12376$ terms. There is a group action of $G \cong S_3 \times S_2 \times S_2$ on $[ABC]$ permuting the $a$, $b$, and $c$ columns appropriately. This extends to a transitive action on the coordinates of the Newton polytope. A few oracle calls quickly determine that $\textrm{New}(X)$ is contained in a $7$-dimensional subspace of $\mathbb{R}^{12}$ and only has $4$ vertices and $2$ facets up to the $G$-action. In total, $\textrm{New}(X)$ has $60$ vertices and $6$ interior points. With only $66$ possible terms, interpolation recovers the polynomial found in  Proposition $7.5$ of \cite{Ponce2017}.
\section{Acknowledgements} We want to express our gratitute to Frank Sottile and Jonathan Hauenstein for many illuminating conversations during the course of this project. We also want to thank Bernd Sturmfels for suggesting the example in Section \ref{sec:algVision} and finally Yue Ren for his suggestion of applying our software to tropical geometry as well as his help with the tropical geometry literature. This work was supported by NSF grant DMS-1501370 and completed during the ICERM-2018 semester on nonlinear algebra.

{
\bibliographystyle{spmpsci.bst}
\bibliography{Bibliography}
}
\end{document}